\documentclass[a4paper,fleqn]{cas-sc}

\usepackage{amsmath}
\usepackage{dsfont}

\usepackage[dvipsnames,svgnames]{xcolor}
\usepackage{hyperref}
\hypersetup{
    colorlinks=true,
    citecolor=RoyalBlue,
    linkcolor=orange,
    filecolor=magenta,      
    urlcolor=cyan
    }

\usepackage[authoryear,longnamesfirst]{natbib}

\usepackage[english]{babel}
\newtheorem{theorem}{Theorem}
\newtheorem{proposition}[theorem]{Proposition}

\newtheorem{corollary}[theorem]{Corollary}

\def\tsc#1{\csdef{#1}{\textsc{\lowercase{#1}}\xspace}}
\tsc{WGM}
\tsc{QE}
\tsc{EP}
\tsc{PMS}
\tsc{BEC}
\tsc{DE}

\begin{document}
\let\WriteBookmarks\relax
\def\floatpagepagefraction{1}
\def\textpagefraction{.001}
\shorttitle{A note on $k$-NN gating in RAG}
\shortauthors{G.~Biau and C.Boyer}

\title [mode = title]{A note on $k$-NN gating in RAG}

\author[1,3]{G\'{e}rard Biau}[orcid=0000-0001-8238-4471]
\ead{gerard.biau@sorbonne-universite.fr}

\affiliation[1]{organization={Sorbonne Université, CNRS, LPSM, F-75005 Paris},
                country={France}}

\author[2,3]{Claire Boyer}[orcid=0000-0002-1566-3371]
\cormark[1]

\ead{claire.boyer@universite-paris-saclay.fr}

\affiliation[2]{organization={Université Paris-Saclay, CNRS, LMO, 91400 Orsay},
country={France}}

\affiliation[3]{organization={Institut universitaire de France},
country={France}}

\cortext[cor1]{Corresponding author}

\begin{abstract}
We propose a statistical proxy framework for retrieval-augmented generation (RAG) that formalizes how language models balance internal predictions with retrieved evidence. We derive an optimal query-level gate, analyze hallucination via retrieval discordance, model query-memory mismatch, and validate the framework numerically on synthetic and real data.
\end{abstract}


\maketitle

\section{Introduction}
Modern language models (LMs) can hallucinate, producing outputs that sound convincing but are factually wrong \citep{ji2023survey,Kalai2024}. Retrieval-augmented generation (RAG) mitigates this by enriching predictions with information drawn from an external memory of documents, code repositories, or previously answered queries \citep{lewis2020rag}. At inference time, the system retrieves items close to the query in an embedding space and conditions its response on them. While effective in improving accuracy and grounding, this mechanism raises a central design question: for each query, how much should the system rely on the LM versus the retrieved evidence?

Most current RAG systems address this balance through heuristics. A common approach concatenates the top-$k$ retrieved items into the prompt \citep{lewis2020rag}, while others interpolate model and retrieval outputs using a fixed mixture weight, as in cache-based or $k$NN language models \citep{grave2017unbounded,khandelwal2020knn,Xu23}. Although often effective, these strategies lack adaptive control: when retrieved neighbors are noisy or off-topic, excessive reliance on retrieval may degrade accuracy, while dominant reliance on the LM can lead to underuse of factual evidence and increased hallucination risk. A principled, query-dependent mechanism for balancing the two sources of information is therefore needed.

This note develops a simple and mathematically explicit framework for studying this balancing problem. The system consists of a frozen base predictor representing the LM, a $k$-nearest neighbor retriever built from an external memory of labeled examples, and a gate that mixes the two. To make the gating decision data-driven and interpretable, we introduce a 
\emph{retrieval-trust weight} $w_{\mathrm{fact}}$ that quantifies how well 
the retrieved neighborhood geometrically supports the query; it acts as a 
penalty in training, discouraging retrieval when local evidence is unreliable. 
This weight, combined with the LM's disagreement with retrieved evidence, 
defines a hallucination discordance criterion \citep{ji2023survey,Kalai2024} 
that the optimal gate is shown to control.

The model we study is not meant to mirror the architectural details of modern RAG systems---prompt construction, attention over retrieved documents, and so on---but to provide an analytically tractable proxy that captures the essential statistical forces governing retrieval-based correction. Numerical experiments in Section~\ref{sec:numerics} and Appendix~\ref{app:experiments} illustrate the resulting gating behavior on synthetic and real data.

\section{Model setup}
We consider an input-output pair \((X,Y)\) drawn from an unknown distribution on 
\(\mathbb R^d\times\mathscr Y\), where \(\mathscr Y=\{1,\ldots,C\}\) is a finite label set. 
An external memory $
\mathscr M_n=\{(U_i,V_i)\}_{i=1}^n$ with 
$(U_i,V_i)\sim Q_{UV}$
stores i.i.d.~reference pairs with \(U\in\mathbb R^d\) and \(V\in\mathscr Y\).  
Throughout the analysis we may assume that the memory is drawn from the same
distribution as \((X,Y)\)---the \emph{aligned setting}, in which \(P_X=Q_U\) and
\(p_{Y|X}=q_{V|U}\).
In practice, however, the memory may come from a related but
distinct source, such as previously answered queries or labeled documents. 
For simplicity, we focus on the aligned setting first and return to this more general situation in Section~\ref{sec:misaligned}.
Retrieval operates in the feature space \(\mathbb R^d\) equipped with a norm $\|\cdot\|$.
For any query \(x\in\mathbb R^d\), let 
\(U_{(1)}(x),\ldots,U_{(k)}(x)\) denote its \(k\) nearest neighbors among 
\(\{U_i\}_{i=1}^n\), with associated labels \(V_{(1)}(x),\ldots,V_{(k)}(x)\). (Ties are broken deterministically by index order.) 

\medskip \noindent\textbf{Base predictor and retriever distribution.}
The base language model (LM) is frozen throughout and provides a conditional 
probability distribution \(p_{\rm LM}(\cdot\mid x)\) on \(\mathscr Y\), interpreted as an 
estimate of the law of \(Y\) given \(X=x\). 
On the retrieval side, the external memory induces its own conditional structure:  
for any query \(x\), the retriever constructs a local, nonparametric estimate of the 
distribution of \(V\) given \(U=x\) by averaging the labels of the $k$ nearest neighbors 
of $x$ in the memory \citep{gyorfi2006distribution,biaudevroye2015}:
\[
\hat r^{(k)}_y(x)
=\frac{1}{k}\sum_{j=1}^k \mathds 1_{\{V_{(j)}(x)=y\}},
\qquad y\in\mathscr Y.
\]
We write $\hat r^{(k)}(y\mid x)=\hat r^{(k)}_y(x)$ for $y\in\mathscr Y$, 
and refer to $\hat r^{(k)}(\cdot\mid x)$ as the \emph{retriever distribution}, 
which targets $p_{Y\mid X}(\cdot\mid x)$ in the aligned setting  
(as made precise by Proposition~\ref{prop:mode-stability} below).

\medskip \noindent\textbf{Retrieval-trust weight.}
To quantify how well the retrieved neighborhood reflects the query, 
we define the retrieval-trust weight
\begin{equation}\label{eq:wret}
w_{\mathrm{fact}}(x)
=\frac{1}{k}\sum_{j=1}^k
\exp\!\big(-\|x-U_{(j)}(x)\|^2\big).
\end{equation}
The weight $w_{\mathrm{fact}}(x)$ is close to $1$ when the retrieved neighbors are close to $x$ and decreases toward $0$ as they spread; it provides a continuous, geometry-aware assessment of retrieval reliability at $x$.

\medskip \noindent\textbf{Mixture model.}
At the core of our framework lies a gated mixture that blends the base LM with the 
retriever. For each query \(x\), predictions are produced by a convex combination of the LM and the retriever distribution:
\begin{equation}\label{eq:mixture}
p_{{\rm RAG}, \lambda}(y\mid x)
=(1-\lambda(x))\,p_{\rm LM}(y\mid x)+\lambda(x)\,\hat r^{(k)}_y(x),
\qquad y\in\mathscr Y,
\end{equation}
where the (measurable) gate \(\lambda:\mathbb R^d\to[0,1]\) controls the 
relative reliance on retrieval.  
Small values of \(\lambda(x)\) favor the LM (fluency and generalization),
while values near one defer to retrieved evidence (grounding and factuality).
Intermediate values achieve an adaptive trade-off.

\medskip \noindent\textbf{Population objective.}
The population-level loss balances predictive accuracy with trust-dependent regularization:
\begin{equation}\label{eq:objective-ret}
\mathscr L(\lambda)
=
\mathbb E \Big[
\sum_{y\in\mathscr Y}p_{Y\mid X}(y\mid X)\,(-\log p_{{\rm RAG}, \lambda}(y\mid X))
\Big]
+\zeta\,\mathbb E \Big[\lambda(X)\,(1-w_{\mathrm{fact}}(X))\Big],
\end{equation}
where $p_{Y\mid X}$ is the true conditional distribution of $Y|X$ and $\zeta\geqslant 0$. The first term enforces predictive fit; the second penalizes retrieval in regions of weak geometric support, with $\zeta$ controlling the penalty strength.

\section{Per-query optimization and hard gating}
The population loss~\eqref{eq:objective-ret} can be written as $\mathscr L(\lambda)=\mathbb E[J(\lambda;X)]$, where
\[
J(\lambda;x)=\sum_{y\in\mathscr Y}p_{Y\mid X}(y\mid x)\,(-\log p_{{\rm RAG}, \lambda}(y\mid x))+\zeta\,\lambda(x)\,(1-w_{\mathrm{fact}}(x)).
\]
Since $\lambda$ enters $J$ only through its value at $x$, minimizing $\mathscr L$ reduces to minimizing $J(\lambda;x)$ pointwise. In the simplest interpretable case, $\lambda(x)\in\{0,1\}$: the mixture~\eqref{eq:mixture} then selects either $p_{\rm LM}$ or $\hat r^{(k)}$, and $J(\lambda;x)$ becomes a two-point comparison between the LM and retriever costs. Defining the local cross-entropies
\[
\ell_{\rm LM}(x)=\sum_{y\in\mathscr Y}p_{Y\mid X}(y\mid x)\,(-\log p_{\rm LM}(y\mid x)),\quad
\ell_r(x)=\sum_{y\in\mathscr Y}p_{Y\mid X}(y\mid x)\,(-\log \hat r^{(k)}_y(x)),
\]
the cost at $\lambda(x)=0$ is $\ell_{\rm LM}(x)$ and the cost at $\lambda(x)=1$ is $\ell_r(x)+\zeta(1-w_{\mathrm{fact}}(x))$.

\begin{proposition}[Optimal hard gate]
\label{prop:hardgate}
For each query \(x\in\mathbb R^d\), the Bayes-optimal hard gate is
\[
\lambda^\star(x)
=
\begin{cases}
0, & \text{if } \quad \ell_{\rm LM}(x)\leqslant \ell_r(x)+\zeta\,(1-w_{\mathrm{fact}}(x)),\\
1, & \text{if } \quad \ell_r(x)+\zeta\,(1-w_{\mathrm{fact}}(x))<\ell_{\rm LM}(x).
\end{cases}
\]
\end{proposition}

The rule of Proposition~\ref{prop:hardgate} states that retrieval is selected exactly when its cross-entropy improvement over the LM exceeds the geometric penalty \(\zeta(1-w_{\mathrm{fact}}(x))\): the decision jointly reflects model fit and neighborhood quality. The parameter $\zeta$ acts as a global regularizer---large values suppress retrieval and favor the LM, small values permit more aggressive grounding in memory.

\medskip \noindent\textbf{Soft gating.}
A continuous gate \(\lambda(x)\in[0,1]\) is also possible: since \(-\log p_{{\rm RAG}, \lambda}(y\mid x)\) is convex in \(\lambda\), the per-query objective \(J(\lambda;x)\) admits a unique minimizer that can be obtained by a one-dimensional numerical solve. We focus on the hard gate of Proposition~\ref{prop:hardgate} because it captures the essential structure of the gating mechanism and facilitates the theoretical analysis below.

\medskip \noindent\textbf{Practical choice of hyperparameters.}
The gate \(\lambda\) is not a free parameter: it is determined by Proposition~\ref{prop:hardgate} once \(\zeta\) is fixed. The neighborhood size \(k\) can be chosen by the standard \(k\)-NN rule \(k=\lfloor\sqrt{n}\rfloor\), which ensures \(k\to\infty\) and \(k/n\to0\) as \(n\to\infty\) \citep{gyorfi2006distribution,biaudevroye2015}, while the regularization \(\zeta\) is naturally selected by cross-validation on a held-out portion of the memory. We illustrate this protocol on real data in Appendix~\ref{app:experiments}.

\section{Hallucination and discordance analysis}
\label{sec:hallucination}
Hallucination---output not supported by the retrieved evidence---is operationalized 
through two ingredients: a measure of how reliable the retrieved evidence is at the 
query, and a measure of how much the LM disagrees with that evidence. 
The retrieval-trust weight $w_{\mathrm{fact}}(x)$ from~\eqref{eq:wret} captures 
the first ingredient (weight of evidence)---large when neighbors are geometrically 
close to $x$, small otherwise---and the quantity $1-p_{\rm LM}(y_r(x)\mid x)$ the second 
(LM's disagreement with retrieved evidence), where $y_r(x)\in{\arg\max}_{y\in\mathscr Y}\ 
\hat r^{(k)}_y(x)$ is the retriever's modal label (ties broken deterministically).

We accordingly define the local discordance score
\(
\mathscr H_{\mathrm{disc}}(p_{\rm LM};x)
=
w_{\mathrm{fact}}(x)\,(1-p_{\rm LM}(y_r(x)\mid x)).
\)
This criterion is large precisely when the retrieval neighborhood is geometrically reliable ($w_{\mathrm{fact}}(x)\approx 1$) but the LM assigns low probability to the label favored by retrieval. In this regime, LM predictions conflict with locally supported evidence, which we interpret as a risk of hallucination.

\subsection{Change under optimal gating}
Under the mixture model~\eqref{eq:mixture},
the hallucination score becomes
\(
\mathscr H_{\mathrm{disc}}(p_{{\rm RAG}, \lambda};x)
= w_{\mathrm{fact}}(x)\,(1-p_{{\rm RAG}, \lambda}(y_r(x)\mid x)).
\)
Thus, the variation relative to the frozen LM is
\begin{equation}
\Delta\mathscr H(x;\lambda) =
\mathscr H_{\mathrm{disc}}(p_{\rm LM};x)-\mathscr H_{\mathrm{disc}}(p_{{\rm RAG}, \lambda};x) = \lambda(x)\,w_{\mathrm{fact}}(x)\,
(\hat r^{(k)}_{y_r(x)}(x)-p_{\rm LM}(y_r(x)\mid x)),
\label{eq:deltaH}
\end{equation}
which is linear in $\lambda(x)$ and satisfies
\(
|\Delta\mathscr H(x;\lambda)|\leqslant \lambda(x)\,w_{\mathrm{fact}}(x).
\)
So, the sign of $\Delta\mathscr H(x;\lambda)$ indicates whether gating 
reduces ($\geqslant 0$) or increases ($\leqslant 0$) local discordance.

Recall that, under hard gating, the optimal decision rule from Proposition~\ref{prop:hardgate} is
\begin{equation}\label{eq:hard-gate-disc}
\lambda^\star(x)
=\mathds 1_{\{\ell_r(x)+\zeta(1-w_{\mathrm{fact}}(x))<\ell_{\rm LM}(x)\}},
\end{equation}
where $\ell_{\rm LM}(x)$ and $\ell_r(x)$
are the LM and retriever local cross-entropies, respectively.
Substituting~\eqref{eq:hard-gate-disc} into~\eqref{eq:deltaH} yields the realized
pointwise change
\begin{equation}\label{eq:deltaH-star}
\Delta\mathscr H(x;\lambda^\star)
=
\mathds 1_{\{\ell_r(x)+\zeta(1-w_{\mathrm{fact}}(x))<\ell_{\rm LM}(x)\}}\,
w_{\mathrm{fact}}(x)\,(\hat r^{(k)}_{y_r(x)}(x)-p_{\rm LM}(y_r(x)\mid x)).
\end{equation}

\medskip \noindent\textbf{Interpretation via three regimes.}
Equation~\eqref{eq:deltaH-star} reveals three qualitatively distinct behaviors governing how gating affects hallucination.

\medskip
\noindent\emph{(i) Gain region.}
On
\(
\mathscr A
=\{\ell_r+\zeta(1-w_{\mathrm{fact}})<\ell_{\rm LM}, \, \
\hat r^{(k)}_{y_r}\geqslant p_{\rm LM}(y_r)\},
\)
the retriever achieves a lower cross-entropy while assigning at least as much mass to its own top label as the LM, so $\Delta\mathscr H(x;\lambda^\star)\geqslant 0$, with magnitude bounded above by $\mathscr H_{\mathrm{disc}}(p_{\rm LM};x)$ and largest when $w_{\mathrm{fact}}(x)\approx 1$.

\medskip
\noindent\emph{(ii) Trade-off region.}
On
\(
\mathscr B
=\{\ell_r+\zeta(1-w_{\mathrm{fact}})<\ell_{\rm LM},\
\hat r^{(k)}_{y_r}<p_{\rm LM}(y_r)\},
\)
the retriever lowers cross-entropy but places less mass on its modal label than the LM, so $\Delta\mathscr H(x;\lambda^\star)\leqslant 0$: a tension between likelihood and factual alignment which the penalty $\zeta(1-w_{\mathrm{fact}}(x))$ mitigates when retrieval is geometrically unreliable.

\medskip
\noindent\emph{(iii) No-switch region.}
On
\(
\mathscr C
=\{\ell_r+\zeta(1-w_{\mathrm{fact}})\geqslant\ell_{\rm LM}\},
\)
the LM remains active ($\lambda^\star(x)=0$) and $\Delta\mathscr H(x;\lambda^\star)=0$.  
This region corresponds to sparse or out-of-distribution queries where retrieval cannot overcome its geometric penalty.

\subsection{Asymptotic analysis of discordance}
The pointwise change $\Delta\mathscr H(x;\lambda^\star)$ in~\eqref{eq:deltaH-star} is governed by the inner quantity $\Delta(x)=\hat r^{(k)}_{y_r(x)}(x)-p_{\rm LM}(y_r(x)\mid x)$, whose sign determines whether gating reduces or increases the local hallucination score. On $\mathscr A$, $\Delta(x)\geqslant 0$ is a clear gain; on $\mathscr B$, $\Delta(x)<0$ despite a cross-entropy improvement, and it is a priori unclear whether this reflects a structural mismatch between the Bayes predictor and the LM, or mere finite-sample variability of the $k$-NN estimator. We address this through a finite-sample mode stability result, from which the asymptotic behavior of $\Delta\mathscr H(x;\lambda^\star)$ follows.

For notational convenience, let $C:=|\mathscr Y|$.  
For $x\in\mathbb R^d$ and $k\ge1$, we denote by
\(
R_k(x)\ =\ \max_{1\leqslant j\leqslant k}\|U_{(j)}(x)-x\|
\)
the $k$-nearest-neighbor radius of the query $x$ among the database points.
\begin{proposition}[Finite-sample mode stability]
\label{prop:mode-stability}
Fix $x\in\mathrm{supp}(Q_U)$ and assume that the conditional distribution $p_{Y|X}(\cdot\mid\cdot)$ 
is $L$-Lipschitz in its second argument. Then, for all $\delta\in(0,1)$,
\[
\mathbb P\Big(
\max_{y\in\mathscr Y}
\big|
\hat r^{(k)}_y(x)
- p_{Y|X}(y\mid x)
\big|
>\delta
\Big)
\leqslant\;  
2C\,\exp \Big(
-2k\Big(\tfrac{\delta}{2}\Big)^{\!2}
\Big)
\;+\;
\mathbb P\Big(R_k(x)>\tfrac{\delta}{2L}\Big).
\]
In particular, if $k\to\infty$ and $k/n \to 0$ as $n\to\infty$, then
\(
\max_{y\in\mathscr Y}
\big|
\hat r^{(k)}_y(x)
- p_{Y|X}(y\mid x)
\big|
\;\stackrel{\mathbb P}{\longrightarrow}\; 0.
\)
\end{proposition}

This proposition provides a uniform finite-sample bound on the deviation $\hat r^{(k)}(\cdot\mid x)-p_{Y|X}(\cdot\mid x)$ at a fixed query~$x$, valid whenever $p_{Y|X}$ is locally
Lipschitz and $x$ lies in the support of the retrieval distribution. In particular, the proposition
shows that, as soon as $k$ grows while $k/n\to 0$, the $k$-NN estimate concentrates around its
Bayes target uniformly over labels. This uniform control is precisely what is needed to guarantee
that, for large samples, the empirical ordering of the coordinates of $\hat r^{(k)}(\cdot\mid x)$ matches the
ordering of the Bayes vector $p_{Y|X}(\cdot\mid x)$. The next corollary makes this consequence
explicit by showing that the empirical top label $y_r(x)$ converges in probability to the Bayes-optimal
label $y^{\star}(x)$ whenever the latter is unique.
\begin{corollary}[Asymptotic mode stability]
\label{cor:mode-stability-consistency}
Fix $x\in\mathrm{supp}(Q_U)$. Assume that the Bayes label
$y^{\star}(x)$ is unique, so that
\[
\gamma(x)
\ :=\
\max_{y\in\mathscr Y}p_{Y|X}(y\mid x)
\;-\;
\max_{y\neq y^{\star}(x)}p_{Y|X}(y\mid x)
\ >\ 0.
\]
Then, under the conditions of Proposition~\ref{prop:mode-stability}, if $k\to\infty$ 
and $k/n\to0$ as $n\to\infty$, 
\(
\mathbb P(y_r(x)\neq y^{\star}(x))\;\longrightarrow\;0.
\)
\end{corollary}

\medskip \noindent\textbf{Asymptotic behavior of the local discordance $\Delta(x)$.}
Combining Proposition~\ref{prop:mode-stability} and Corollary~\ref{cor:mode-stability-consistency}, we analyze the large-sample behavior of the key quantity $\Delta(x)$, which determines the sign and magnitude of the hallucination variation $\Delta\mathscr H(x;\lambda^\star)$ in~\eqref{eq:deltaH-star}. Fix \(x\in\mathrm{supp}(Q_U)\) and assume that the Bayes label \(y^{\star}(x)\) is unique, so that
\(\gamma(x)>0\). By Proposition~\ref{prop:mode-stability}, if \(k\to\infty\) and \(k/n\to0\), then
\(
\max_{y\in\mathscr Y}\big|\hat r^{(k)}_y(x)-p_{Y|X}(y\mid x)\big|
\;\stackrel{\mathbb P}{\longrightarrow}\;0.
\)
Therefore
\(
\hat r^{(k)}_{y_r(x)}(x)-p_{Y|X}(y_r(x)\mid x)
\;\stackrel{\mathbb P}{\longrightarrow}\;0.
\)
Recalling \(\Delta(x)=\hat r^{(k)}_{y_r(x)}(x)-p_{\rm LM}(y_r(x)\mid x)\), this implies $
\Delta(x)
-
\big(p_{Y|X}(y_r(x)\mid x)-p_{\rm LM}(y_r(x)\mid x)\big)
\;\stackrel{\mathbb P}{\longrightarrow}\;0$.
Moreover, Corollary~\ref{cor:mode-stability-consistency} yields $
\mathbb P\big(y_r(x)=y^{\star}(x)\big)\;\longrightarrow\;1$. On the event \(\{y_r(x)=y^{\star}(x)\}\),
\[
p_{Y|X}(y_r(x)\mid x)-p_{\rm LM}(y_r(x)\mid x)
=
p_{Y|X}(y^{\star}(x)\mid x)-p_{\rm LM}(y^{\star}(x)\mid x).
\]
Hence,
\[
p_{Y|X}(y_r(x)\mid x)-p_{\rm LM}(y_r(x)\mid x)
\;\stackrel{\mathbb P}{\longrightarrow}\;
p_{Y|X}(y^{\star}(x)\mid x)-p_{\rm LM}(y^{\star}(x)\mid x),
\]
and combining with the previous display gives
\begin{equation}
\label{asymp-delta}
\Delta(x)
\;\stackrel{\mathbb P}{\longrightarrow}\;
p_{Y|X}(y^{\star}(x)\mid x)-p_{\rm LM}(y^{\star}(x)\mid x).
\end{equation}
In particular, the sign of \(\Delta(x)\) coincides with that of
\(p_{Y|X}(y^{\star}(x)\mid x)-p_{\rm LM}(y^{\star}(x)\mid x)\) with probability tending to one.

\medskip \noindent\textbf{Asymptotic behavior of the local hallucination variation.}
We now turn to the asymptotic behavior of the local hallucination variation $\Delta\mathscr H(x;\lambda^\star)$ in~\eqref{eq:deltaH-star}.
The result below shows that, for fixed $x$ in the support,
the sign and magnitude of $\Delta\mathscr H(x;\lambda^\star)$ converge in probability
to a deterministic quantity governed solely by the Bayes predictor and the LM.
We let $\ell_{\mathrm{Bayes}}(x)=\sum_{y\in\mathscr Y}p_{Y|X}(y\mid x)\,(-\log p_{Y|X}(y\mid x))$ and recall that the LM local cross-entropy is $\ell_{\rm LM}(x)=\sum_{y\in\mathscr Y}p_{Y|X}(y\mid x)\,(-\log p_{\rm LM}(y\mid x))$.
\begin{theorem}[Asymptotic behavior of the local hallucination variation]
\label{thm:asymptotic-deltaH}
Assume the aligned setting, and let $x\in\mathrm{supp}(Q_U)$.
Suppose that the Bayes label $
y^{\star}(x)\in\arg\max_{y\in\mathscr Y}p_{Y|X}(y\mid x)$ is unique and that $p_{Y|X}(\cdot\mid\cdot)$
is $L$-Lipschitz in its second argument. 
Then, if $k\to\infty$ and $k/n\to0$ as $n\to\infty$,
\[
\Delta\mathscr H(x;\lambda^\star)
\;\stackrel{\mathbb P}{\longrightarrow}\;
\lambda_\infty(x)\,
(p_{Y|X}(y^{\star}(x)\mid x)-p_{\rm LM}(y^{\star}(x)\mid x)),
\]
where $
\lambda_\infty(x)
\ :=\
\mathds 1_{\{\ell_{\mathrm{Bayes}}(x)<\ell_{\rm LM}(x)\}}$.
\end{theorem}

The variation $\Delta\mathscr H(x;\lambda^\star)$ therefore converges in probability to a quantity determined entirely by the structural relationship between the Bayes rule and the LM at~$x$: the randomness of both the $k$-NN estimator and the trust weight $w_{\mathrm{fact}}(x)$ vanishes in the limit. When $\ell_{\mathrm{Bayes}}(x)<\ell_{\rm LM}(x)$ the gate activates with probability tending to one and $\Delta\mathscr H(x;\lambda^\star)\to p_{Y|X}(y^{\star}(x)\mid x)-p_{\rm LM}(y^{\star}(x)\mid x)$, which may take either sign; when $\ell_{\mathrm{Bayes}}(x)=\ell_{\rm LM}(x)$ the gate remains off and $\Delta\mathscr H(x;\lambda^\star)\to 0$. Any persistent nonzero variation in the large-sample regime is therefore structural rather than a $k$-NN finite-sample artifact.

\section{Numerical illustration}
\label{sec:numerics}

We illustrate the gating analysis on a controlled two-dimensional synthetic setting designed to exhibit the three regimes $\mathscr{A}$, $\mathscr{B}$, $\mathscr{C}$ identified in Section~\ref{sec:hallucination}. We take $\mathscr Y=\{1,2,3\}$, draw the memory $\{(U_i,V_i)\}_{i=1}^{2000}$ with $U_i$ uniform on $[0,1]^2$ and $V_i\sim p_{Y|X}(\cdot\mid U_i)$, where $p_{Y\mid X}(y\mid x)\propto\exp(-(x_1-c_y)^2/\tau)$ with $c_1=0.2$, $c_2=0.5$, $c_3=0.8$, $\tau=0.05$. The base model $p_{\rm LM}$ is designed piecewise in $x_2$ to exhibit each regime: 
on $\{x_2<1/3\}$ we set $p_{\rm LM}=p_{Y\mid X}$ (well-specified LM); 
on $\{1/3\leqslant x_2<2/3\}$ we set $p_{\rm LM}$ to a cyclic permutation of $p_{Y\mid X}$ (wrong modal class); and on $\{x_2\geqslant 2/3\}$ we set $p_{\rm LM}\propto p_{Y\mid X}^{4}$ (overconfident on the right modal class). We take $k=30$ and $\zeta=2$.

\begin{figure}
\centering
\includegraphics[width=\linewidth]{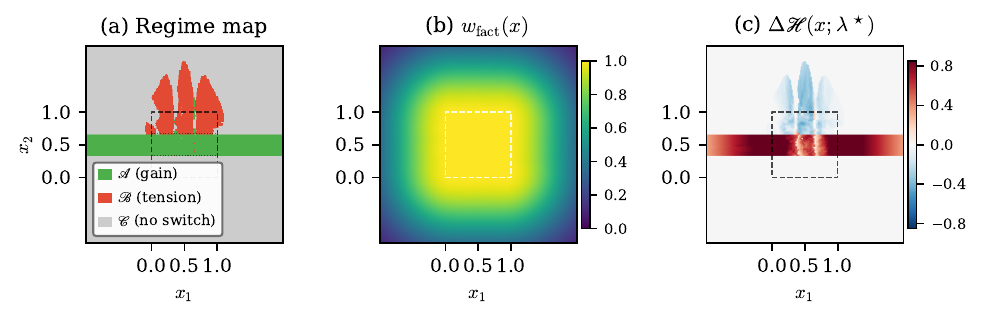}
\caption{Three-regime visualization on the 2D synthetic setting described in Section~\ref{sec:numerics}. (a) Regime map of $\lambda^\star$, with $\mathscr{A}$ (gain, green), $\mathscr{B}$ (tension, red), $\mathscr{C}$ (no switch, grey). (b) Retrieval-trust weight $w_{\mathrm{fact}}(x)$. (c) Local discordance variation $\Delta\mathscr{H}(x;\lambda^\star)$ from~\eqref{eq:deltaH-star}, with positive values in red ($\mathscr{A}$) and negative values in blue ($\mathscr{B}$), white at zero. Memory support $[0,1]^2$ shown dashed.}
\label{fig:regimes}
\end{figure}

Figure~\ref{fig:regimes} confirms the theoretical decomposition. The three $x_2$-indexed horizontal bands visible inside the dashed square $[0,1]^2$ correspond to the three $p_{\rm LM}$ regimes, with the middle band ($\mathscr{A}$) showing $\Delta\mathscr{H}>0$ (gating reduces discordance) and the top band ($\mathscr{B}$) showing $\Delta\mathscr{H}<0$ (the tension case). In the bottom band, $p_{\rm LM}=p_{Y\mid X}$ implies $\ell_r(x)\geqslant\ell_{\rm LM}(x)$ and the gate 
never fires, yielding a fluctuation-free $\mathscr{C}$ region for any~$n$. The $\mathscr{A}$ region extends laterally beyond $[0,1]^2$: the cross-entropy gain from retrieval survives a moderate geometric penalty when the LM is severely misspecified. As $\|x\|$ grows further, $w_{\mathrm{fact}}(x)$ decays and the penalty drives the gate into $\mathscr{C}$, so the geometric component acts as a built-in safeguard against out-of-support retrieval. Further experiments---asymptotic convergence of $\Delta\mathscr{H}$ to the limit of Theorem~\ref{thm:asymptotic-deltaH}, and a real-data study with cross-validated $\zeta$ and a misaligned variant---are reported in Appendix~\ref{app:experiments}.

\section{Beyond the aligned setting}
\label{sec:misaligned}
The analysis above assumed $P_X=Q_U$ and $p_{Y|X}=q_{V|U}$. In practice, queries and memory typically follow distinct laws, with possible covariate and semantic shifts between them \citep{Zhou2023,tamang2025survey}. Appendix~\ref{sec:mismatch} develops a hybrid geometric-semantic mismatch model that captures both effects through a perturbation $X=T(U)=U+\xi(U)$ of memory inputs and a local corruption of memory labels; here we record its central consequence for the trust weight $w_{\mathrm{fact}}$.

Let $S=\operatorname{supp}(Q_U)$ and $d(x,S)=\inf_{u\in S}\|x-u\|$ be the distance from $x$ to the memory support.
\begin{proposition}[Asymptotic behavior of the trust weight]
\label{prop:wfact-limit}
Fix $x\in\mathbb R^d$. If $k/n\to0$ as $n\to\infty$, then
\[
w_{\mathrm{fact}}(x)
=\frac{1}{k}\sum_{j=1}^{k}\exp\!\big(-\|x-U_{(j)}(x)\|^2\big)
\;\stackrel{\mathrm{a.s.}}{\longrightarrow}\;\exp\!\big(-d(x,S)^2\big).
\]
In particular, $w_{\mathrm{fact}}(x)\stackrel{\mathrm{a.s.}}{\to}1$ for $x\in S$.
\end{proposition}

This gives the penalty $\zeta(1-w_{\mathrm{fact}}(x))$ in~\eqref{eq:objective-ret} an unambiguous geometric meaning under any form of mismatch: it grows precisely with how far the query lies from where memory was collected, irrespective of whether the discrepancy is covariate, semantic, or combined. The gate is therefore automatically steered toward the base model whenever retrieval would draw from an unsupported region of feature space. Appendix~\ref{sec:mismatch} formalizes the corresponding limit of the retriever distribution $\hat r^{(k)}$ and quantifies the joint geometric-semantic bias.

\medskip\noindent{\textit{Acknowledgments.} We thank the anonymous referee for constructive feedback that substantially improved the paper.}

\bibliographystyle{plainnat}
\bibliography{biblio}

\begin{thebibliography}{11}
\providecommand{\natexlab}[1]{#1}
\providecommand{\url}[1]{\texttt{#1}}
\expandafter\ifx\csname urlstyle\endcsname\relax
  \providecommand{\doi}[1]{doi: #1}\else
  \providecommand{\doi}{doi: \begingroup \urlstyle{rm}\Url}\fi

\bibitem[Alpaydin and Kaynak(1998)]{alpaydin1998digits}
Ethem Alpaydin and Cenk Kaynak.
\newblock Optical recognition of handwritten digits.
\newblock {UCI} Machine Learning Repository, 1998.

\bibitem[Biau and Devroye(2015)]{biaudevroye2015}
G\'erard Biau and Luc Devroye.
\newblock \emph{Lectures on the Nearest Neighbor Method}.
\newblock Springer, Cham, 2015.

\bibitem[Grave et~al.(2017)Grave, Cisse, and Joulin]{grave2017unbounded}
Edouard Grave, Moustapha~M Cisse, and Armand Joulin.
\newblock Unbounded cache model for online language modeling with open
  vocabulary.
\newblock In Isabelle Guyon, Ulrike von Luxburg, Samy Bengio, Hanna Wallach,
  Rob Fergus, S.~Vishwanathan, and Roman Garnett, editors, \emph{Advances in
  Neural Information Processing Systems}, volume~30, pages 6044--6054. Curran
  Associates, Inc., 2017.

\bibitem[Gy{\"o}rfi et~al.(2006)Gy{\"o}rfi, Kohler, Krzy{\.z}ak, and
  Walk]{gyorfi2006distribution}
L{\'a}szl{\'o} Gy{\"o}rfi, Michael Kohler, Adam Krzy{\.z}ak, and Harro Walk.
\newblock \emph{A Distribution-Free Theory of Nonparametric Regression}.
\newblock Springer, New York, 2006.

\bibitem[Ji et~al.(2023)Ji, Lee, Frieske, Yu, et~al.]{ji2023survey}
Ziwei Ji, Nayeon Lee, Rita Frieske, Tiezheng Yu, et~al.
\newblock Survey of hallucination in natural language generation.
\newblock \emph{ACM Computing Surveys}, 55:\penalty0 248,1--38, 2023.

\bibitem[Kalai and Vempala(2024)]{Kalai2024}
Adam~Tauman Kalai and Santosh~S. Vempala.
\newblock Calibrated language models must hallucinate.
\newblock In \emph{Proceedings of the 56th Annual ACM Symposium on Theory of
  Computing}, STOC 2024, pages 160--171, New York, 2024. Association for
  Computing Machinery.

\bibitem[Khandelwal et~al.(2020)Khandelwal, Levy, Jurafsky, Zettlemoyer, and
  Lewis]{khandelwal2020knn}
Urvashi Khandelwal, Omer Levy, Dan Jurafsky, Luke Zettlemoyer, and Mike Lewis.
\newblock Generalization through memorization: {N}earest neighbor language
  models.
\newblock In \emph{International Conference on Learning Representations}, 2020.

\bibitem[Lewis et~al.(2020)Lewis, Perez, Piktus, Petroni, et~al.]{lewis2020rag}
Patrick Lewis, Ethan Perez, Aleksandra Piktus, Fabio Petroni, et~al.
\newblock Retrieval-augmented generation for knowledge-intensive {NLP} tasks.
\newblock In Hugo Larochelle, Marc{'}Aurelio Ranzato, Raia Hadsell,
  Maria-Florina Balcan, and Haibin Lin, editors, \emph{Advances in Neural
  Information Processing Systems}, volume~33, pages 9459--9474. Curran
  Associates, Inc., 2020.

\bibitem[Tamang et~al.(2025)Tamang, Bouadjenek, Dazeley, and
  Aryal]{tamang2025survey}
Lakpa Tamang, Mohamed~Reda Bouadjenek, Richard Dazeley, and Sunil Aryal.
\newblock Handling out-of-distribution data: {A} survey.
\newblock \emph{arXiv:2507.21160}, 2025.

\bibitem[Xu et~al.(2023)Xu, Alon, and Neubig]{Xu23}
Frank~F. Xu, Uri Alon, and Graham Neubig.
\newblock Why do nearest neighbor language models work?
\newblock In Andreas Krause, Emma Brunskill, Kyunghyun Cho, Barbara Engelhardt,
  Sivan Sabato, and Jonathan Scarlett, editors, \emph{Proceedings of the 40th
  International Conference on Machine Learning}, volume 202 of
  \emph{Proceedings of Machine Learning Research}, pages 38325--38341. JMLR,
  2023.

\bibitem[Zhou et~al.(2023)Zhou, Liu, Qiao, Xiang, and Change]{Zhou2023}
Kaiyang Zhou, Ziwei Liu, Yu~Qiao, Tao Xiang, and Loy~Chen Change.
\newblock Domain generalization: {A} survey.
\newblock \emph{IEEE Transactions on Pattern Analysis and Machine
  Intelligence}, 45:\penalty0 4396--4415, 2023.

\end{thebibliography}

\clearpage

\begin{center}
    \vspace*{0.5cm}
    
    \Large Appendix to "A Note on $k$-NN gating in RAG"
    \vspace*{0.5cm}
\end{center}

\hrulefill

\appendix
\section{Proofs}
\subsection{Proof of Proposition~\ref{prop:mode-stability}}
Conditionally on the neighbor locations $U_{(1)}(x),\ldots,U_{(k)}(x)$, the variables
\[
Z_{j,y}:=\mathds 1_{\{V_{(j)}(x)=y\}},\qquad y\in\mathscr Y,
\]
are independent Bernoulli with means $p_{Y|X}(y\mid U_{(j)}(x))$, and
$\hat r^{(k)}_y(x)=k^{-1}\sum_{j=1}^k Z_{j,y}$.
Therefore, by Hoeffding's inequality and a union bound over all classes,
for any $\delta\in(0,1)$,
\setcounter{equation}{7}
\begin{align}
&\mathbb P \Big(
\max_{y\in\mathscr Y}
\Big|
\hat r^{(k)}_y(x)
- \frac{1}{k}\sum_{j=1}^k p_{Y|X}(y\mid U_{(j)}(x))
\Big|
>\frac{\delta}{2}
\;\Big|\;
U_{(1)}(x),\ldots,U_{(k)}(x)
\Big) \leqslant
2C\,\exp\!\Big(
-2k\Big(\tfrac{\delta}{2}\Big)^2
\Big).
\label{eq:cond-hoeffding-prop}
\end{align}
Dropping the conditioning yields the same bound unconditionally.

Next define the local modulus of continuity
\[
\omega_x(r)
:=\sup_{\|u-x\|\leqslant r}\max_{y\in\mathscr Y}
|p_{Y|X}(y\mid u)-p_{Y|X}(y\mid x)|.
\]
Clearly,
\[
\max_{y\in\mathscr Y}
\Big|
\frac{1}{k}\sum_{j=1}^k p_{Y|X}(y\mid U_{(j)}(x))
- p_{Y|X}(y\mid x)
\Big|
\leqslant \omega_x(R_k(x)).
\]
If $p_{Y|X}$ is $L$-Lipschitz, then $\omega_x(r)\leqslant Lr$, and combining this with  
\eqref{eq:cond-hoeffding-prop} yields
\begin{equation}
\mathbb P\Big(
\max_{y\in\mathscr Y}
\big|
\hat r^{(k)}_y(x)
- p_{Y|X}(y\mid x)
\big|
>\delta
\Big)
\leqslant\;
2C\,\exp\!\Big(
-2k\Big(\tfrac{\delta}{2}\Big)^{\!2}
\Big)
\;+\;
\mathbb P\Big(R_k(x)>\tfrac{\delta}{2L}\Big).
\label{inter}
\end{equation}

Finally, since $x\in\mathrm{supp}(Q_U)$, the local mass 
$q_{x,\varepsilon}:=Q_U(B(x,\varepsilon))$ is strictly positive for every $\varepsilon>0$.
Because $\{R_k(x)>\varepsilon\}=\{\mathrm{Bin}(n,q_{x,\varepsilon})<k\}$, Chernoff's bound \citep[e.g.,][Chapter~20]{biaudevroye2015} gives, whenever  
$k\leqslant \tfrac12\,n q_{x,\varepsilon}$,
\[
\mathbb P\big(R_k(x)>\varepsilon\big)
\ \leqslant\
\exp\!\Big(-\tfrac{n q_{x,\varepsilon}}{8}\Big).
\]
Because $k/n\to0$ implies $k\leqslant \tfrac12\,n q_{x,\varepsilon}$ for large $n$, the right-hand
side tends to zero.  
Combining this with~\eqref{inter} gives the desired convergence.
\subsection{Proof of Corollary~\ref{cor:mode-stability-consistency}}
Fix $x\in\mathrm{supp}(Q_U)$ and assume that the Bayes label $y^{\star}(x)$ is unique, so that $\gamma(x)>0$.
Set $\varepsilon := \tfrac{\gamma(x)}{3}$.
By Proposition~\ref{prop:mode-stability}, applied with this choice of $\varepsilon$, we obtain
\[
\mathbb P\Big(
\max_{y\in\mathscr Y}
\big|
\hat r^{(k)}_y(x)
- p_{Y|X}(y\mid x)
\big|
>\varepsilon
\Big)
\;\longrightarrow\;0
\qquad\text{as }n\to\infty,
\]
whenever $k\to\infty$ and $k/n\to0$. Now, define the event
\[
A_n
\ :=\
\Big\{
\max_{y\in\mathscr Y}
\big|
\hat r^{(k)}_y(x)
- p_{Y|X}(y\mid x)
\big|
\leqslant\varepsilon
\Big\},
\]
so that $\mathbb P(A_n^c)\to 0$ as $n\to\infty$.
On $A_n$, for any $y_1,y_2\in\mathscr Y$,
\[
\hat r^{(k)}_{y_1}(x)-\hat r^{(k)}_{y_2}(x)
\ \geqslant\
p_{Y|X}(y_1\mid x)-p_{Y|X}(y_2\mid x)
\ -\ 2\varepsilon.
\]
In particular, taking $y_1=y^{\star}(x)$ and any $y_2\neq y^{\star}(x)$ gives
\[
p_{Y|X}(y^{\star}(x)\mid x)-p_{Y|X}(y_2\mid x)
\ \geqslant\ \gamma(x)
\ =\ 3\varepsilon,
\]
hence
\[
\hat r^{(k)}_{y^{\star}(x)}(x)-\hat r^{(k)}_{y_2}(x)
\ \geqslant\ 3\varepsilon-2\varepsilon
\ =\ \varepsilon
\ >\ 0.
\]
Thus, on $A_n$,
\[
\hat r^{(k)}_{y^{\star}(x)}(x)
>\hat r^{(k)}_{y}(x)
\qquad\forall\,y\neq y^{\star}(x),
\]
and therefore the empirical and Bayes top labels coincide:
\[
y_r(x)
=\arg\max_{y\in\mathscr Y}\hat r^{(k)}_y(x)
=\arg\max_{y\in\mathscr Y}p_{Y|X}(y\mid x)
=y^{\star}(x).
\]
Consequently,
\[
\mathbb P (y_r(x)\neq y^{\star}(x))
\;\leqslant\;
\mathbb P(A_n^c)
\;\longrightarrow\;0
\qquad\text{as }n\to\infty,
\]
which establishes the claim.
\subsection{Proof of Theorem~\ref{thm:asymptotic-deltaH}}
If $\ell_{\rm LM}(x)=\ell_{\mathrm{Bayes}}(x)$, then KL non-negativity implies
$\ell_r(x)\geqslant\ell_{\mathrm{Bayes}}(x)=\ell_{\rm LM}(x)$ for every memory realization,
so $\lambda^\star(x)=0$ and $\Delta\mathscr H(x;\lambda^\star)=0$ identically:
the conclusion holds trivially.
We therefore assume $\ell_{\rm LM}(x)>\ell_{\mathrm{Bayes}}(x)$,
i.e.\ $\lambda_\infty(x)=1$, throughout.

Since $x\in\mathrm{supp}(Q_U)$ and $k/n\to0$, 
Proposition~\ref{prop:wfact-limit} yields
\begin{equation}
w_{\mathrm{fact}}(x)\stackrel{\mathbb P}{\longrightarrow}1.
\label{asymp-w}
\end{equation}
Next, recall that
\[
\lambda^\star(x)
=\mathds 1_{\{\ell_r(x)+\zeta(1-w_{\mathrm{fact}}(x))<\ell_{\rm LM}(x)\}},
\]
where
\[
\ell_r(x)
=\sum_{y\in\mathscr Y}p_{Y|X}(y\mid x)\,(-\log \hat r^{(k)}_y(x)).
\]
By Proposition~\ref{prop:mode-stability} and the convention $0\log 0=0$,
continuity of $\log$ on labels with $p_{Y|X}(y\mid x)>0$ implies $
\ell_r(x)\;\stackrel{\mathbb P}{\longrightarrow}\;
\ell_{\mathrm{Bayes}}(x)$. Using \eqref{asymp-w}, we obtain
\[
\ell_r(x)+\zeta(1-w_{\mathrm{fact}}(x))
\;\stackrel{\mathbb P}{\longrightarrow}\;
\ell_{\mathrm{Bayes}}(x).
\]

Since \(\ell_{\mathrm{Bayes}}(x)<\ell_{\rm LM}(x)\) in the present case,  we have
\[
\mathds 1_{\{\ell_r(x)+\zeta(1-w_{\mathrm{fact}}(x))<\ell_{\rm LM}(x)\}}
\;\stackrel{\mathbb P}{\longrightarrow}\;1.
\]
Equivalently,
\[
\mathds 1_{\{\ell_r(x)+\zeta(1-w_{\mathrm{fact}}(x))<\ell_{\rm LM}(x)\}}
\;\stackrel{\mathbb P}{\longrightarrow}\;
\lambda_\infty(x),
\qquad
\lambda_\infty(x)
=
\mathds 1_{\{\ell_{\mathrm{Bayes}}(x)<\ell_{\rm LM}(x)\}}.
\]
Combining this limit with~\eqref{asymp-delta} proves the theorem.


\section{Query-memory distribution mismatch: Model and retriever asymptotics}
\label{sec:mismatch}

This appendix formalizes the hybrid mismatch model alluded to in Section~\ref{sec:misaligned} and completes the asymptotic analysis with the limit of the $k$-NN retriever distribution. The aligned setting of the main text corresponds to $\xi\equiv0$ and $\rho\equiv0$ in the model below. Sources of mismatch (covariate and semantic shift) and their statistical study are reviewed in \citet{Zhou2023} and \citet{tamang2025survey}.

\subsection{A hybrid geometric-semantic mismatch model}
We characterize query-memory mismatch by distinguishing the distribution \(P_X\) of
query inputs  from that of memory inputs \(Q_U\), and the query labeling
mechanism \(p_{Y|X}\) from the memory labeling mechanism \(q_{V|U}\).
This perspective allows us to model both geometric distortion in the embedding
space and semantic mismatch in the conditional labels.

\medskip \noindent\textbf{Geometric shift.}
Queries are assumed to be generated from memory inputs through a perturbation model:
\[
X=T(U)=U+\xi(U),
\]
where \(U\sim Q_U\) and \(\xi:\mathbb R^d\to\mathbb R^d\) is a deformation function. The distribution \(P_X=T_\# Q_U\) thus represents the law of queries obtained
by displacing the memory inputs.  
When computing the $k$ nearest neighbors of a query \(x\),
the retrieved points \(U_{(1)}(x),\ldots,U_{(k)}(x)\) are the elements of
\(\{U_i\}_{i=1}^n\) that lie closest to \(x\) in the embedding space.
If $\|\xi(u)\|$ is small, the neighborhoods of $x$ and $u$ largely overlap;
as $\|\xi(u)\|$ increases, retrieved neighbors become less representative of \(x\),
capturing the geometric component of the mismatch.

\medskip \noindent\textbf{Label drift.}
Even when geometric distortion is negligible, the conditional relationship between features 
and labels in memory may differ from that of current queries.
We model this semantic deviation as a local corruption process:
\[
q_{V|U}(y\mid u)
=(1-\rho(u))\,p_{Y|X}(y\mid u)
+\rho(u)\,s(y\mid u),
\qquad 0\leqslant\rho(u)\leqslant1,
\]
where $\rho(u)$ is a local corruption rate and \(s(\cdot\mid u)\) an arbitrary spurious distribution.
Thus the retrieval distribution constructed from \(\mathscr M_n\) approximates a corrupted 
version of \(p_{Y|X}\): the memory is reliable when $\rho(u)$ is small and increasingly 
misleading as $\rho(u)$ grows.

\medskip \noindent\textbf{Retrieval trust and interpretation.}
The two mechanisms combine into the coupled model
\setcounter{equation}{10}
\begin{equation}\label{eq:hybrid-model}
U\sim Q_U,\quad X=T(U),\quad
V\mid U\sim q_{V|U}(\cdot\mid U),\quad
Y\mid X\sim p_{Y|X}(\cdot\mid X).
\end{equation}
Both types of shift influence the reliability of retrieval, but in different ways.
The retrieval-trust weight $w_{\mathrm{fact}}(x)$, defined in~\eqref{eq:wret}, 
measures the geometric compatibility between the query $x$ and its retrieved neighbors.
Large geometric deformations $\|\xi(u)\|$ inflate the distances $\|x-U_{(j)}(x)\|$ and
therefore directly reduce $w_{\mathrm{fact}}(x)$.
By contrast, strong semantic corruption $\rho(u)$ does not affect $w_{\mathrm{fact}}(x)$
itself, but makes the retrieved labels unreliable as proxies for the true query labels,
even when geometric proximity is high.
Thus $w_{\mathrm{fact}}(x)$ should be interpreted as a geometry-aware indicator of
retrieval quality, while the retrieval distribution captures the semantic reliability
of the retrieved labels under joint geometric and semantic shift.

The following subsection formalizes the asymptotic behavior of the $k$-NN 
retriever $\hat r^{(k)}(\cdot\mid x)$ under mild regularity assumptions on 
$T$ and $\rho$. Proposition~\ref{prop:rhat-limit} establishes that the 
retriever converges to the local memory label distribution at the nearest 
support points. Together with Proposition~\ref{prop:wfact-limit} of 
Section~\ref{sec:misaligned}, this gives a precise statistical interpretation 
of the penalty term in~\eqref{eq:objective-ret}: retrieval is discouraged 
precisely in regions where geometric proximity is weak or where the limiting 
retrieval distribution reflects boundary or corrupted semantics.

\subsection{Asymptotic behavior of the retriever under the hybrid model}
Let \(\{(U_i,V_i)\}_{i=1}^n\) be the memory pairs with \(U_i\in\mathbb R^d\), drawn i.i.d.\ from \(Q_{UV}\).
Denote by \(Q_U\) the marginal distribution of the \(U_i\)'s, and by \(S=\operatorname{supp}(Q_U)\) its (closed) support.
For any query \(x\in\mathbb R^d\), we let $ d(x,S)=\inf_{u\in S}\|x-u\|$. (Since $S$ is closed, this infimum is attained.)

Under the geometric component of the hybrid model,
\(X=T(U)=U+\xi(U)\) with \(U\sim Q_U\),
a query \(x=T(u)\) satisfies \(d(x,S)\leqslant \|x-u\|=\|\xi(u)\|\).
Therefore, Proposition~\ref{prop:wfact-limit} yields the almost-sure limit
\[
w_{\mathrm{fact}}(x)
\;\xrightarrow[n\to\infty]{\mathrm{a.s.}}\;
\exp\!\big(-d(x,S)^2\big)
\;\geqslant\;
\exp\!\big(-\|\xi(u)\|^2\big),
\]
with equality whenever \(\|x-u\|=d(x,S)\).
\begin{proposition}[Asymptotics of the $k$-NN retrieval distribution]
\label{prop:rhat-limit}
Fix \(x\in\mathbb R^d\). Assume that \(k\to\infty\) and \(k/n\to 0\) as \(n\to\infty\), and define
\[
\hat r^{(k)}_y(x)=\frac{1}{k}\sum_{j=1}^k \mathds 1_{\{V_{(j)}(x)=y\}},
\qquad y\in\mathscr Y.
\]
Let \(S=\operatorname{supp}(Q_U)\) and let 
\[
N_S(x)=\{u\in S:\|x-u\|=d(x,S)\}
\]
be the (possibly non-singleton) set of nearest points in \(S\) to \(x\).

\smallskip
\noindent\emph{(i) {\bf In-support point.}}
If \(x\in S\) and \(q_{V|U}(y\mid\cdot)\) is continuous at \(x\), then
\[
\hat r^{(k)}_y(x)\; \stackrel{\mathrm{a.s.}}{\longrightarrow} \; q_{V|U}(y\mid x).
\]

\smallskip
\noindent\emph{(ii) {\bf Unique nearest point off support.}}
If \(x\notin S\), the nearest set is a singleton \(N_S(x)=\{u_x\}\), and 
\(q_{V|U}(y\mid\cdot)\) is continuous at \(u_x\), then
\[
\hat r^{(k)}_y(x)\; \stackrel{\mathrm{a.s.}}{\longrightarrow} \; q_{V|U}(y\mid u_x).
\]

\smallskip
\noindent\emph{(iii) {\bf Multiple nearest points.}}
If \(q_{V|U}(y\mid\cdot)\) is continuous on \(N_S(x)\), then almost surely,
\[
\min_{u\in N_S(x)} q_{V|U}(y\mid u) 
\leqslant
\liminf_{n\to\infty}\hat r^{(k)}_y(x)
\leqslant
\limsup_{n\to\infty}\hat r^{(k)}_y(x)
\leqslant
\max_{u\in N_S(x)} q_{V|U}(y\mid u).
\]
In particular, if \(q_{V|U}(y\mid u)\) is constant on \(N_S(x)\), then 
\(\hat r^{(k)}_y(x)\) converges to that common value.
\end{proposition}

\medskip \noindent\textbf{Interpretation under the hybrid shift model.}
Proposition~\ref{prop:rhat-limit} shows that
the empirical retriever \(\hat r^{(k)}(\cdot\mid x)\) consistently estimates the
memory's local label mechanism in the region of the database that is
geometrically closest to the query.  
When \(x\in S\), the estimate converges to \(q_{V|U}(\cdot\mid x)\); when
\(x\notin S\) but admits a unique projection \(u_x\in N_S(x)\), it converges
to \(q_{V|U}(\cdot\mid u_x)\).
Proposition~\ref{prop:wfact-limit} complements this semantic characterization
with a geometric one: \(w_{\mathrm{fact}}(x)\approx 1\) indicates that \(x\) lies in
or close to \(S\), whereas small values of \(w_{\mathrm{fact}}(x)\) flag
out-of-support queries for which retrieval reflects boundary behavior rather
than genuine local structure.
Thus, under geometric and semantic mismatch, the pair
\((\hat r^{(k)}, w_{\mathrm{fact}})\) jointly encodes 
the reliability of retrieved evidence, a property that directly supports and
stabilizes the gating rule.

Propositions~\ref{prop:wfact-limit} and~\ref{prop:rhat-limit} yield almost-sure
pointwise limits for \(w_{\mathrm{fact}}(x)\) and---when the limit exists---for
\(\hat r^{(k)}(\cdot\mid x)\).  In particular, in cases~(i) and~(ii) of
Proposition~\ref{prop:rhat-limit}, the $k$-NN retriever admits a deterministic
limit \(r_\infty(\cdot\mid x)\). Whenever \(r_\infty(\cdot\mid x)\) exists, continuity of the mixture implies that,
for any measurable \(\lambda:\mathbb R^d\to[0,1]\),
\begin{align*}
p_{{\rm RAG}, \lambda}(y\mid x)
&=(1-\lambda(x))\,p_{\rm LM}(y\mid x)
   +\lambda(x)\,\hat r^{(k)}_y(x)\\
&\xrightarrow[n\to\infty]{\mathrm{a.s.}}\ 
p_{\lambda,\infty}(y\mid x)
:=(1-\lambda(x))\,p_{\rm LM}(y\mid x)
  +\lambda(x)\,r_\infty(y\mid x).
\end{align*}
Thus, under the hybrid shift model~\eqref{eq:hybrid-model}, we obtain
\begin{align*}
&r_\infty(y\mid x)-p_{Y|X}(y\mid x)\\
&=
(1-\rho(u_x))\big(p_{Y|X}(y\mid u_x)-p_{Y|X}(y\mid x)\big)
+\rho(u_x)\big(s(y\mid u_x)-p_{Y|X}(y\mid x)\big),
\end{align*}
and, writing \(\|\cdot\|_1=\sum_{y\in\mathscr Y}|\cdot|\),
\[
\|r_\infty(\cdot\mid x)-p_{Y|X}(\cdot\mid x)\|_1
\ \leqslant\
(1-\rho(u_x))\,\delta_{\mathrm{geom}}(x)
+\rho(u_x)\,\delta_{\mathrm{sem}}(x),
\]
where
\[
\delta_{\mathrm{geom}}(x)=\|p_{Y|X}(\cdot\mid u_x)-p_{Y|X}(\cdot\mid x)\|_1,
\quad
\delta_{\mathrm{sem}}(x)=\|s(\cdot\mid u_x)-p_{Y|X}(\cdot\mid x)\|_1.
\]
If \(p_{Y|X}(\cdot\mid x)\) is locally Lipschitz as a map from \(\mathbb R^d\) to
\((\mathbb R^C,\|\cdot\|_1)\), then
\[
\delta_{\mathrm{geom}}(x)
=\|p_{Y|X}(\cdot\mid u_x)-p_{Y|X}(\cdot\mid x)\|_1
\;\leqslant\;
L\,\|x-u_x\|
\;\leqslant\;
L\,d(x,S).
\]
Hence the total retrieval bias is jointly driven by the geometric displacement
\(d(x,S)\) and the local corruption level \(\rho(u_x)\).
When \(x\) lies near the memory support and corruption is small, the limiting
retriever \(r_\infty(\cdot\mid x)\) is close to the true conditional distribution,
and \(w_{\mathrm{fact}}(x)\approx 1\) makes the penalty
\(\lambda(x)(1-w_{\mathrm{fact}}(x))\) negligible, so retrieval is not discouraged.
As \(x\) moves away from the support or corruption increases,
\(\|r_\infty(\cdot\mid x)-p_{Y|X}(\cdot\mid x)\|_1\) grows while
\(w_{\mathrm{fact}}(x)\approx e^{-d(x,S)^2}\) shrinks, naturally steering the
gate toward the base model.

\subsection{Proof of Proposition~\ref{prop:wfact-limit}}
Let \(D_n^{(j)}(x)=\|x-U_{(j)}(x)\|\) be the $j$-th nearest-neighbor distance among 
\(\{U_i\}_{i=1}^n\).  
By Lemma~2.2 of \citet{biaudevroye2015}, if \(k/n\to0\) then
\(D_n^{(k)}(x) \stackrel{\mathrm{a.s.}}{\longrightarrow} d(x,S)\) as $n\to \infty$. Since \(D_n^{(1)}(x)\leqslant\cdots\leqslant D_n^{(k)}(x)\) and \(D_n^{(1)}(x)\geqslant d(x,S)\), we obtain
\[
0\leqslant
\max_{1\leqslant j\leqslant k}\big|D_n^{(j)}(x)-d(x,S)\big|
\leqslant
D_n^{(k)}(x)-d(x,S)
\,\stackrel{\mathrm{a.s.}}{\longrightarrow} \, 0,
\]
so \(D_n^{(j)}(x)\stackrel{\mathrm{a.s.}}{\longrightarrow} d(x,S)\) uniformly over \(1\leqslant j\leqslant k\).
With \(\varphi(t)=\exp(-t^2)\) continuous, uniform convergence gives
\[
\max_{1\leqslant j\leqslant k}\big|\varphi(D_n^{(j)}(x))-\varphi(d(x,S))\big|
\,\stackrel{\mathrm{a.s.}}{\longrightarrow} \, 0.
\]
This shows the desired claim.

\subsection{Proof of Proposition~\ref{prop:rhat-limit}}
Let \(D_n^{(j)}(x)=\|x-U_{(j)}(x)\|\).  As in the proof of Proposition~\ref{prop:wfact-limit}, when \(k/n\to 0\) as $n\to \infty$,
\[
0\ \leqslant\ 
\max_{1\leqslant j\leqslant k}|D_n^{(j)}(x)-d(x,S)|
\ \leqslant\ 
D_n^{(k)}(x)-d(x,S)
\,\stackrel{\mathrm{a.s.}}{\longrightarrow} \, 0.
\]
Thus the neighbor locations satisfy almost surely: if \(x\in S\), then \(U_{(j)}(x)\to x\); if \(x\notin S\) and \(N_S(x)=\{u_x\}\), then \(U_{(j)}(x)\to u_x\);  in general, every cluster point of the sequence \(\{U_{(j)}(x)\}_{j=1}^k\) lies in \(N_S(x)\). 

Conditionally on the neighbor locations, the variables  
\(Z_{j,y}:=\mathds 1_{\{V_{(j)}(x)=y\}}\) are independent Bernoulli with means
\(q_{V|U}(y\mid U_{(j)}(x))\), and
$\hat r^{(k)}_y(x)=k^{-1}\sum_{j=1}^k Z_{j,y}$.
Thus, Hoeffding's inequality yields
\[
\mathbb P \Big(
\Big|\hat r^{(k)}_y(x) - \frac{1}{k}\sum_{j=1}^k q_{V|U}(y\mid U_{(j)}(x))\Big|
>\varepsilon
\;\Big|\;\ U_{(1)}(x),\ldots,U_{(k)}(x)
\Big)
\leqslant 2e^{-2k\varepsilon^2},
\]
hence, as $k\to \infty$,
\[
\hat r^{(k)}_y(x) - \frac{1}{k}\sum_{j=1}^k q_{V|U}(y\mid U_{(j)}(x))
\,\stackrel{\mathrm{a.s.}}{\longrightarrow} \, 0.
\]

For (i), continuity at \(x\in S\) implies
\[
\max_{1\leqslant j\leqslant k}
|q_{V|U}(y\mid U_{(j)}(x)) - q_{V|U}(y\mid x)|
\,\stackrel{\mathrm{a.s.}}{\longrightarrow} \, 0,
\]
so the Cesàro mean converges to \(q_{V|U}(y\mid x)\).  
The same argument gives (ii).

For (iii), continuity on \(N_S(x)\) implies the Cesàro averages of  
\(\{q_{V|U}(y\mid U_{(j)}(x))\}_{j=1}^k\) must lie within the convex hull of the
function values on \(N_S(x)\), giving the stated bounds.  

\section{Numerical experiments}
\label{app:experiments}

This appendix complements the synthetic illustration of Section~\ref{sec:numerics} with two further experiments: a real-data study on the UCI handwritten digits dataset that demonstrates cross-validated selection of $\zeta$ and robustness to covariate shift, and a finite-sample validation of the asymptotic limit of Theorem~\ref{thm:asymptotic-deltaH}.

\subsection{Real data: {C}ross-validated gating on digits}
\medskip \noindent\textbf{Setup.} 
The UCI handwritten digits dataset \citep{alpaydin1998digits} comprises 1797 grayscale images of 10 digit classes, each summarized by 64 pixel features which we $\ell_2$-normalize. We split the 1797 samples into a memory $\mathscr M_n=\{(U_i,V_i)\}_{i=1}^{1257}$ 
(70\%) and a held-out test set. The memory is further partitioned into a retrieval memory $\mathscr M_{\mathrm{retr}}$ (80\%, 1006 samples) and a validation set $\mathscr M_{\mathrm{val}}$ (20\%, 251 samples) used to select $\zeta$. The neighborhood size is set to $k=\lfloor\sqrt{|\mathscr M_{\mathrm{retr}}|}\rfloor=31$, in line with the standard $k$-NN rule.
The base model $p_{\rm LM}$ is a multinomial logistic regression trained on $\mathscr M_{\mathrm{retr}}$ 
with $30\%$ of the labels replaced by uniformly random ones, mimicking a deliberately 
miscalibrated LM. The hyperparameter $\zeta$ is then selected by maximizing accuracy on $\mathscr M_{\mathrm{val}}$ over a logarithmic grid in $[10^{-2}, 10^{1.5}]$. Note that, for simplicity, $p_{\rm LM}$ and $\hat r^{(k)}$ are built from the same pool
$\mathscr M_{\mathrm{retr}}$; the label corruption models a LM that has seen
the domain but with imperfect supervision, while the retriever holds the
curated labels---the intended RAG scenario.

\begin{figure}
\centering
\includegraphics[width=\linewidth]{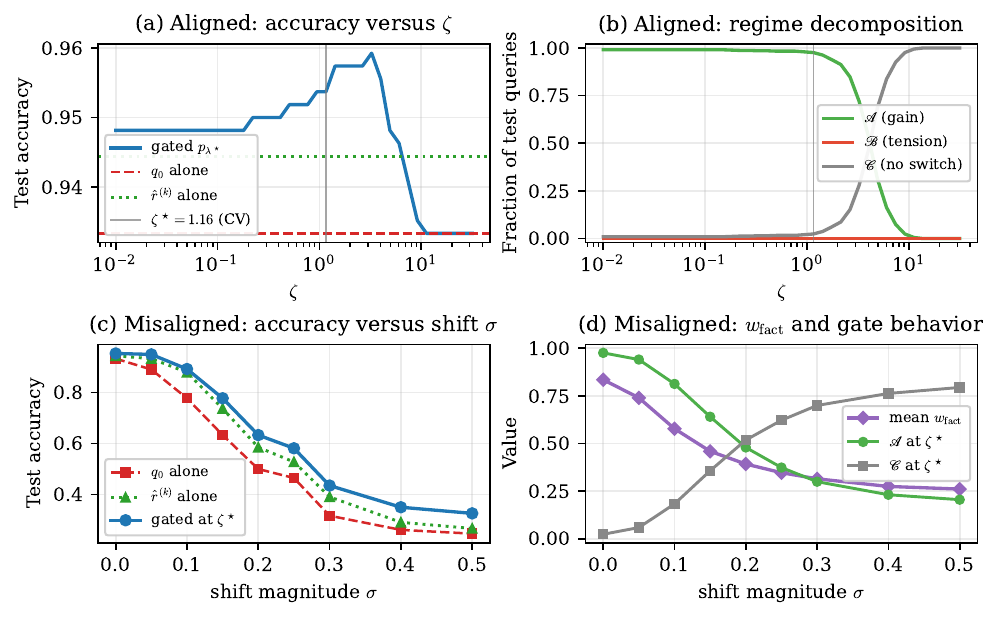}
\caption{Real-data experiment on UCI digits. (a) Test accuracy of the gated mixture as $\zeta$ varies, compared with $p_{\rm LM}$ alone (dashed) and $\hat r^{(k)}$ alone (dotted); the CV-selected $\zeta^\star=1.16$ is marked. (b) Regime decomposition on the test set as $\zeta$ varies. (c) Test accuracy under covariate shift, where each test query is perturbed by additive Gaussian pixel noise of standard deviation $\sigma$, the gate using the same $\zeta^\star$. (d) Mean retrieval-trust weight $w_{\mathrm{fact}}$ and gate composition at $\zeta^\star$ as $\sigma$ grows.}
\label{fig:realdata}
\end{figure}

\medskip \noindent\textbf{Results.} Figure~\ref{fig:realdata}(a) shows that the gated mixture $p_{\lambda^\star}$ strictly outperforms both baselines across a broad range of $\zeta$. With $\zeta^\star=1.16$ selected by cross-validation, the test accuracy of the mixture is $0.954$, against $0.944$ for $\hat r^{(k)}$ alone and $0.933$ for $p_{\rm LM}$ alone. Panel (b) shows the corresponding regime decomposition: at $\zeta^\star$ the gate fires on a vast majority of queries (regime $\mathscr A$), reflecting the substantial miscalibration of $p_{\rm LM}$. Regime $\mathscr{B}$ is negligible throughout; it is theoretically meaningful 
but its conditions are rarely met in this experiment.

Panels (c) and (d) examine robustness to distribution shift: each test query is 
perturbed by additive Gaussian noise of standard deviation $\sigma$ followed by 
$\ell_2$-renormalization, displacing it on the feature sphere while its class label 
remains that of the original image. The gate uses the same $\zeta^\star$ selected 
on the aligned validation set. As $\sigma$ grows, both baselines degrade rapidly while the gated mixture remains uniformly above them; the gap reaches roughly $5$ percentage points above $\hat r^{(k)}$ alone at $\sigma=0.3$, where $\hat r^{(k)}$ achieves $0.39$, $p_{\rm LM}$ achieves $0.32$, and the gate achieves $0.44$. Panel (d) tracks the mean $w_{\mathrm{fact}}$ and the fractions $|\mathscr A|/n_{\mathrm{test}}$ and $|\mathscr C|/n_{\mathrm{test}}$ at $\zeta^\star$: $w_{\mathrm{fact}}$ decreases from $0.83$ to $0.26$ as $\sigma$ increases from $0$ to $0.5$, and the gate gracefully shifts mass from $\mathscr A$ to $\mathscr C$. This is the empirical counterpart of Proposition~\ref{prop:wfact-limit}: $w_{\mathrm{fact}}$ tracks the distance of the (perturbed) query to the memory support, and the geometric penalty automatically steers the gate toward the base model when retrieval becomes unreliable.

\subsection{Synthetic validation of Theorem~\ref{thm:asymptotic-deltaH}}
We finally validate the asymptotic statement of Theorem~\ref{thm:asymptotic-deltaH} on the synthetic setup of Section~\ref{sec:numerics}. For each $n\in\{50,100,200,500,1000,2000,5000,10000\}$ we set $k=\lfloor\sqrt n\rfloor$, draw $M=200$ independent memory realizations, and compute $\Delta\mathscr{H}(x;\lambda^\star)$ at three representative query points, one per regime: $x_{\mathscr A}=(0.40,0.50)$ (label-permuted LM), $x_{\mathscr B}=(0.40,0.85)$ (sharpened LM), $x_{\mathscr C}=(0.40,0.15)$ (well-specified LM). The theoretical limits given by Theorem~\ref{thm:asymptotic-deltaH} evaluate to $+0.282$, $-0.291$, and $0$, respectively.

\begin{figure}
\centering
\includegraphics[width=\linewidth]{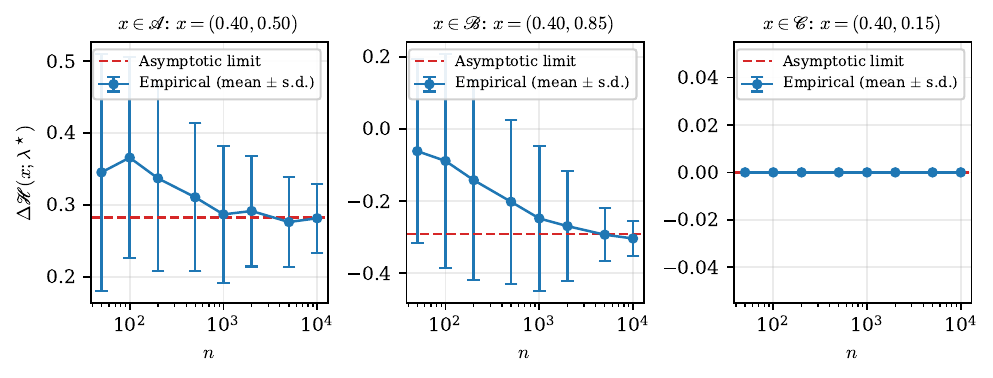}
\caption{Empirical mean and standard deviation of $\Delta\mathscr{H}(x;\lambda^\star)$ across $M=200$ independent memory realizations, for three query points representative of each regime (left to right: $\mathscr A$, $\mathscr B$, $\mathscr C$). Asymptotic limits from Theorem~\ref{thm:asymptotic-deltaH} are shown dashed.}
\label{fig:scaling}
\end{figure}

Figure~\ref{fig:scaling} reports the empirical mean and standard deviation of $\Delta\mathscr{H}(x;\lambda^\star)$ as a function of $n$, with the predicted limit overlaid. In all three regimes the empirical mean converges to the theoretical limit, and the standard deviation shrinks at the expected rate. In regime $\mathscr C$, $p_{\rm LM}=p_{Y\mid X}$ implies $\lambda_\infty(x)=0$;
KL non-negativity further forces $\ell_r(x)\geqslant\ell_{\rm LM}(x)$ for every
memory realization, so the gate never fires and $\Delta\mathscr H(x;\lambda^\star)=0$
exactly for all~$n$.

\end{document}